\documentclass[letterpaper, 11pt,  reqno]{amsart}

\usepackage{amsmath,amssymb,amscd,amsthm,amsxtra, esint}

\usepackage[implicit=true]{hyperref}

\usepackage{color}
\usepackage{ulem}
\usepackage[makeroom]{cancel}

\allowdisplaybreaks[2]

\newtheorem{theorem}{Theorem} [section]

\newtheorem{remark}[theorem]{Remark}

\newtheorem{definition}{Definition}

\DeclareMathOperator*{\supp}{supp}

\newcommand{\nn}{|\hspace{-0.4mm}|\hspace{-0.4mm}|}

\newcommand{\noi}{\noindent}
\newcommand{\Z}{\mathbb{Z}}
\newcommand{\R}{\mathbb{R}}

\newcommand{\F}{\mathcal{F}}

\newcommand{\al}{\alpha}

\newcommand{\dl}{\delta}

\newcommand{\nb}{\nabla}

\newcommand{\eps}{\varepsilon}

\newcommand{\ld}{\lambda}

\newcommand{\s}{\sigma}

\newcommand{\ft}{\widehat}

\newcommand{\cj}{\overline}

\renewcommand{\l}{\ell}

\newcommand{\les}{\lesssim}
\newcommand{\ges}{\gtrsim}

\newcommand{\jb}[1]
{\langle #1 \rangle}

\newcommand{\ind}{\mathbf 1}

\renewcommand{\S}{\mathcal{S}}

\newcommand{\M}{\mathfrak{M}}

\newcommand{\ww}{\mathfrak{w}}

\newcommand{\NB}{\mathbb{N}}

\newtheorem*{ackno}{Acknowledgements}

\numberwithin{equation}{section}
\numberwithin{theorem}{section}

\newcommand{\w}{\pmb{w}}

\renewcommand{\ss}{\pmb{\sigma}}

\begin{document}
\baselineskip = 14pt

\title{Modulation spaces with scaling symmetry}
\author{\'Arp\'ad B\'enyi \and Tadahiro Oh}

\address{\'Arp\'ad B\'enyi, Department of Mathematics,
Western Washington University, 516 High Street, Bellingham, WA 98225,
USA}
\email{arpad.benyi@wwu.edu}

\address{Tadahiro Oh, School of Mathematics,
The University of Edinburgh and The Maxwell Institute for the Mathematical Sciences,
James Clerk Maxwell Building, The King's Buildings, Peter Guthrie Tait Road,
Edinburgh, EH9 3FD, UK}

\email{hiro.oh@ed.ac.uk}




\begin{abstract}
We indicate how to construct a family of modulation spaces that have a scaling symmetry. We also illustrate the behavior of the Schr\"odinger multiplier on such function spaces.
\end{abstract}

\maketitle

\section{Reconciling the modulation and dilation scalings}

The Besov spaces and the modulation spaces are both examples of so-called decomposition spaces that stem from either a dyadic covering or uniform covering of the underlying frequency space. Roughly speaking, both of these classes of function spaces are defined by imposing appropriate decay conditions on a sequence of the form
$\big\{\|\mathcal F^{-1}(\psi_k\widehat{f})\|_{L^p_x}\big\}_{k}$,
where $\{\psi_k\}_{k\in \Z}$ is a dyadic partition of unity in the Besov case,
while $\{\psi_k\}_{k\in\Z^d}$ is a uniform partition of unity in the case of the modulation spaces. Moreover, one can connect the two classes of spaces via the so-called $\alpha$-modulation spaces of Gr\"obner \cite{G} for $\alpha\in [0, 1]$; the modulation
spaces (and the Besov spaces, respectively) are obtained at the end-point: $\alpha=0$
(and $\alpha=1$,
respectively).
The $\alpha$-modulation spaces are still decomposition spaces obtained via a ``geometric interpolation" method that now asks for a covering $\{Q_k^\alpha\}_{k \in I}$ of the frequency space in which the sets satisfy a condition of the form $|Q_k^\al|\sim \langle \xi\rangle^{\alpha}$ for all $\xi\in Q_k^\alpha$; see \cite{FV}.

In this paper,  however, we are interested in a different connection between the two classes of function spaces than that alluded to above. It is clear from their definitions that the Besov spaces enjoy the dilation symmetry: $f(x) \mapsto f(2^k x)$, $k \in \Z$,
while the modulation spaces enjoy the modulation symmetry:
$f(x) \mapsto e^{2\pi i k  \cdot x} f(x)$, $k \in \Z^d$.
A natural question then is if there exists a class of function spaces that is sufficiently rich that would enjoy both such symmetries. 
As we shall see, reconciling the modulation and dilation symmetries is possible by considering scaled versions of the modulation spaces and then appropriately amalgamating them with a ``good" vector weight.
 In many respects, such a class of function spaces would be more enticing to consider from the perspective of someone working in PDEs. Our note stems from 
 these natural considerations and should be seen as a contribution to the general program of constructing new function spaces from given ones, as well as deriving the essential properties of the new spaces 
 (i) from those of the modulation spaces on which they are rooted and 
 (ii) from the characteristics of the weight, as proposed in \cite{Fei2015}.

\section{Definition of the $\M^{p, q, r}_{\w}$-spaces and their scaling property}

\subsection{The $\M^{p, q, r}_{\w}$-spaces}
Let $p, q, r\in [1, \infty)$ be fixed. We denote by $\w$ a vector weight $\{w_j\}_{j\in\Z}$ with $w_j\geq0$.
Generally speaking, we aim to define a space which captures the modulations of all scales, while the choice of the weight ${\w}$ is made according to the task for which the space is needed.

Let us first recall next the definition of modulation spaces $M^{p, q}$; see \cite{FG1, FG2}. Let $\psi \in \S(\R^d)$ such that
\begin{align}
\supp \psi \subset Q_0:=[-1, 1]^d
\qquad \text{and} \qquad \sum_{k \in \Z^d} \psi(\xi -k) \equiv 1.
\label{psi0}
\end{align}

\noi
Then, the modulation space $M^{p, q}$ is defined
as the collection of all tempered distributions
$f\in\S'(\R^d)$ such that
$\|f\|_{M^{p, q}}<\infty$, where
the $M^{p, q}$-norm is defined by
\begin{equation} \label{mod}
\|f\|_{M^{p, q} (\R^d)} = \big\|\|\psi(D-k) f
\|_{L_x^p(\R^d)} \big\|_{\l^q_k(\Z^d)}.
\end{equation}

\noi
Here, $\psi (D-k)f(x)=\int_{\R^d}\psi (\xi-k)\ft f(\xi)e^{2\pi ix\cdot \xi}\,d\xi$. Clearly, the support of $\psi$ renders the domain of the integration of the previous integral to be $\xi\in Q_0+k$, and in particular we have
\begin{align}
 f = \sum_{k\in\Z^d}\psi(D-k)f.
 \label{decomp1}
\end{align}

In what follows, for ease of notation and unless specifically stated otherwise, we assume that the underlying space is $\R^d$. Now, for a (fixed) \textit{dyadic scale} $j\in \mathbb Z$, 
we define the scaled modulation space $M^{p,q}_{[j]}$ by the norm:
\[\|f\|_{M^{p,q}_{[j]}}:=\big\|\|\psi_{j, k}(D) f
\|_{L_x^p(\R^d)} \big\|_{\l^q_k(\Z^d)},\]

\noi
where $\psi_{j, k}(D)$ is defined by
\begin{align}
\psi_{j, k}(D)f(x)=\int_{\R^d}\psi (2^{-j}\xi-k)\ft f(\xi)e^{2\pi ix\cdot \xi}\,d\xi.
\label{psi1}
\end{align}

\noi
Namely, while the usual modulation spaces $M^{p, q}$ are adapted
to the modulation symmetry of the unit scale,
the $M^{p,q}_{[j]}$-spaces are adapted to the modulation symmetry of scale $2^j$.
When $j = 0$, we have  $M^{p, q}_{[0]}=M^{p,q}$.
Note also that, as in \eqref{decomp1}, we have
\begin{align}
f = \sum_{k\in\Z^d} \psi_{j, k}(D)f
\label{decomp2}
\end{align}

\noi
for all $j\in \Z$.

Now, in view of the decompositions \eqref{decomp1} and \eqref{decomp2},
we have
\begin{align*}
f = \sum_{j \in \Z}  \sum_{k\in\Z^d}w_j \psi_{j, k}(D)f
\end{align*}

\noi
for any vector weight $\w = \{w_j \}_{j\in \Z}$ with $\|\w\|_{\l^1}=1$.
This motivates the following definition of the modulation spaces with scaling symmetry.

\begin{definition}\label{DEF:mod_sc}\rm
Let $1\leq p, q, r<\infty$
and  $\w = \{w_j\}_{j\in\Z}$ be a vector weight with $w_j\geq0$.
We define the space $\M^{p, q, r}_{\w}(\R^d)$
as the collection of all tempered distributions $f\in \S'(\R^d)$ such that $\|f\|_{\M^{p,q, r}_{\w}}<\infty$, where
the $\M^{p,q, r}_{\w}$-norm is defined by
\begin{equation}\label{mod-new}
\|f\|_{\M^{p,q, r}_{\w}}:=\big\|\|f\|_{M^{p,q}_{[j]}}\big\|_{\ell^r_j(\w)}
=\bigg(\sum_{j\in\Z}
w_j^r\|f\|_{M^{p,q}_{[j]}}^r\bigg)^{\frac 1r}.
\end{equation}
\end{definition}

Note that the definition of the space
$\M^{p,q, r}_{\w}$ heavily depends on the choice of the weight~$\w$
and we need to make sure that smooth functions belong to
$\M^{p,q, r}_{\w}$ for suitable weights $\w$.
We first claim  that the Schwartz class $\S(\R^d)$ is contained in $\M^{p, q, r}_{\w}$,
provided that  $\w=\{w_j\}_{j\in\Z}$ is a ``good'' vector weight, that is, for some $\eps>0$ we have
\begin{align}
w_j\les
\begin{cases}
2^{-\eps j} & \text{if}\,\,j\geq 0,\\
2^{-   (\frac{d}{p'} - \frac dq - \eps) j } & \text{if}\,\,j<0.
\end{cases}
\label{X1}
\end{align}

\noi
Note that we allow some growth for $j < 0$, when $q > p'$.
For $j \geq 0$, it suffices to impose
summability of $\{w_j \}_{j \geq 0}$.

Let us first prove this claim. Let $\psi_{j, k}$ denote
the multiplier of $\psi_{j, k}(D)$ defined in \eqref{psi1},
corresponding to the smoothed version of $\ind_{_{2^j(Q_0+k)}}$.
Now, for $j<0$, it follows from Young's inequality that
\begin{align*}
\|\psi_{j, k}(D)f\|_{L^p_x}
& \leq  \|\mathcal F^{-1}(\psi_{j, k})\|_{L^p_x} \bigg\|\sum_{|i| \le1 }  \psi_{j, k+i} (D) f\, \bigg\|_{L^1_x} \notag\\
& \sim 2^{j \frac{d}{p'}} \bigg\|\sum_{|i| \le1 }  \psi_{j, k+i} (D) f\, \bigg\|_{L^1_x},
\end{align*}

\noi
where the implicit constant is independent of $k \in \Z^d$.
By computing the $\l^q_k$-norm, we have
\begin{align}
\|f\|_{M^{p,q}_{[j]}}
& = \big\| \|\psi_{j, k}(D)f\|_{L^p_x}\big\|_{\l^q_k(\Z^d)}
 \les 2^{j \frac{d}{p'}} \big\|\|  \psi_{j, k} (D) f\|_{L^1_x} \big\|_{\l^q_k(\Z^d)}\notag\\
&  \leq 2^{j \frac{d}{p'}}
 \big\| \jb{2^j k}^{-s}\| \psi_{j, k} (D)  \jb{\nb}^s  f \|_{L^1_x}\big\|_{\l^q_k(\Z^d)}\notag\\
 &  \les
 2^{j \frac{d}{p'}} \|f\|_{W^{s, 1}_x}
 \bigg( \sum_{k \in \Z^d} \frac{1}{\jb{2^j k}^{sq}}\bigg)^\frac{1}{q}\notag\\
&    \les
 2^{j d (\frac{1}{p'} - \frac 1q)} \|f\|_{W^{s, 1}_x}
\label{X3}
\end{align}

\noi
for any $s > \frac dq$.
Here, $\jb{x}  = (1 + |x|^2)^\frac{1}{2}$
and we used the Riemann sum approximation in the last step.
Then,
from \eqref{X1} and \eqref{X3},
the contribution to
the $\M^{p,q, r}_{\w}$-norm for $j < 0$ is estimated by
\begin{align*}
\bigg(\sum_{j < 0}w_j^r
\|f\|_{M^{p,q}_{[j]}}^r
\bigg)^\frac{1}{r}
\les
\Big(\sum_{j < 0}2^{\eps j r} \Big)^\frac{1}{r}
 \|f \|_{W^{s, 1}_x}
<\infty.
\end{align*}

When $j \geq 0$, proceeding as above,
we have
\begin{align}
\|f\|_{M^{p,q}_{[j]}}
    \les
 \|f\|_{W^{s, p}_x}
\label{X4}
\end{align}

\noi
for any $s > \frac dq$.
Hence, from \eqref{X1} and \eqref{X4}, we obtain
\begin{align*}
\bigg(\sum_{j \geq 0}w_j^r
\|f\|_{M^{p,q}_{[j]}}^r
\bigg)^\frac{1}{r}
\les
\Big(\sum_{j \geq  0}2^{- \eps j r} \Big)^\frac{1}{r}
 \|f \|_{W^{s, p}_x}
<\infty.
\end{align*}

We point out that  the condition \eqref{X1} is  sharp.
Let $j < 0$.
Suppose that there exists $J \in \NB$ such that
 $w_j\gtrsim 2^{-(\frac d{p'} - \frac dq) j}$ for any $j\leq -J$.
Let $f \in \S(\R^d)$.
Then, there exist  $ N \in \NB$ and $k_0 \in \Z^d$ such that
\[|\ft f(\xi)| \ge \frac 12 \|\ft f \|_{L^\infty}\]

\noi
for all $\xi \in 2^{-N} (Q_0 + k_0)$
and $\ft f$ essentially behaves like a constant on the cube $2^{-N} (Q_0 + k_0)$.
This in particular implies that
\[
\|\psi_{j, k}(D)f\|_{L^p_x}\ges\|\ft f\|_{L^\infty}2^{j\frac d{p'}},
\]

\noi
for any $j \leq -N$ and $k \in \Z^d$ such that
\begin{align}
2^j (Q_0 + k ) \subset 2^{-N} (Q_0 + k_0).
\label{X5}
\end{align}

\noi
Then, it is easy to see that  $\|f\|_{\M^{p, q, r}_{\w}}=\infty$
by simply computing the contribution from $j < 0$ and $k \in \Z^d$ satisfying \eqref{X5}.
The condition \eqref{X1}
for $j \geq 0$ can also be seen to be essentially sharp by noting that
\[ \|f\|_{M^{p,q}_{[j]}}
\geq \|\psi_{j, 0}(D) f\|_{L_x^p} \ges  \| f\|_{L^p_x}\]

\noi
for all sufficiently large $j \geq 0$.

Lastly, we wish to summarize the essential properties of the function spaces $\M^{p, q, r}_{\w}$ which bear a strong resemblance with those of the classical modulation spaces. The dilation property will be treated separately in the next subsection. See also Section \ref{SEC:L^2} for a connection to other spaces stemming from PDEs.

\begin{theorem}
\label{properties}
Let  $1\leq p, q, r<\infty$ and assume that the vector weight $\w$ satisfies \eqref{X1}. 
Then, we have the following:

\noi 
\begin{itemize}
\item
[\textup{(i)}] $\S(\R^d)\subset \M^{p, q, r}_{\w}(\R^d)$, 

\item[\textup{(ii)}] $\M^{p, q, r}_{\w}(\R^d)\subset M^{p, q}(\R^d)=M^{p, q}_{[0]}(\R^d)$,

\item[\textup{(iii)}] $\big(\M^{p,q, r}_{\w}(\R^d), \|\cdot\|_{\M^{p,q, r}_{\w}}\big)$ is a Banach space, 
\rule[-2.5mm]{0pt}{0pt}

\item[\textup{(iv)}] $\M^{p_1, q_1, r_1}_{\w}(\R^d)\subset \M^{p_2, q_2, r_2}_{\w}(\R^d)$ for $1\leq p_1\leq p_2$, $1\leq q_1\leq q_2$, and $1\leq r_1\leq r_2$.
\end{itemize}
\end{theorem}

\begin{remark}\rm

A natural and interesting question arises, concerning the duality of our function spaces. 
Let $1 \leq  p, q, r < \infty$ and 
fix  a vector weight $\w = \w(p, q) = \{w_j \}_{j \in \Z}$ satisfying the condition \eqref{X1}.
 Then, 
by setting 
\begin{align*}
 \jb{f, g}
: = \sum_{j \in \Z} w_j w'_j \sum_{k \in \Z^d} \int_{\R^d} \big(\psi_{j, k}(D)f\big)(x)  \big(\psi_{j, k}(D)g\big)(x)dx
\end{align*}

\noi
for $f \in \M^{p, q, r}_{\w}$, 
H\"older's inequality yields
\begin{align*}
| \jb{f, g}|
\leq \| f\|_{\M^{p, q, r}_{\w}}
\|g\|_{\M^{p', q', r'}_{\w'}}, 
\end{align*}

\noi
where  $\w' = \w'(p', q') = \{w'_j\}_{j \in \Z}$ is a vector weight satisfying the condition \eqref{X1} 
with $(p', q')$ in place of $(p, q)$.
We call such a vector weight $\w'$ a dual weight.
This shows that 
$\bigcup_{\w'} \M^{p', q', r'}_{\w'} \subset \big(\M^{p, q, r}_{\w}\big)'$,
where the union is taken over all the dual weights $\w'$ (for a given pair $(p, q)$).
It is, however, not clear to us how  to determine
the actual dual space $\big(\M^{p, q, r}_{\w}\big)'$ in this formulation.

Let us discuss an alternative approach.
Let $m$ be a counting measure on $\Z$
and write \eqref{mod-new}
as
\[ 
\|f\|_{\M^{p,q, r}_{\w}}:=\| \psi_{j, k}(D) f\|_{L^r (\Z, w_j^r dm)\l^q_k (\Z^d) L^p_x(\R^d)}. 
\]

\noi
Then, a duality pairing may be given by 
\begin{align*}
(f, g)
: = \int  \bigg(\sum_{k \in \Z^d} \int_{\R^d} \big(\psi_{j, k}(D)f\big)(x)  \big(\psi_{j, k}(D)g\big)(x)dx\bigg)
\, w_j^r  dm
\end{align*}

\noi
In this case, H\"older's inequality (viewing $w_j^r  dm$ as a measure on $\Z$) yields
\begin{align*}
|(f, g)| \leq \| f\|_{\M^{p, q, r}_{\w}}
\|g\|_{\M^{p', q', r'}_{\w'}}, 
\end{align*}

\noi
where $\w' = \{ w_j^{r - 1}\}_{j \in \Z}$.
This shows $ \M^{p', q', r'}_{\w'} \subset \big(\M^{p, q, r}_{\w}\big)'$
for this particular weight $\w'$.
Note that this weight $\w'$ may not satisfy the condition \eqref{X1}
unless we impose a further assumption on $\w$.

\end{remark}

\subsection{Scaling  property}
\label{SUBSEC:scaling}

Let us now investigate the scaling property of
the new modulation spaces $\M^{p,q, r}_{\w}$.
Define a translation operator $\tau$ on a vector weight  $\w = \{ w_j \}_{j \in \Z}$
by setting $(\tau \w)_j = w_{j+1}$.
For $k\in\NB$, write $\tau^{k}$ (and $\tau^{-k}$, respectively) for
$\tau$ composed with itself $k$ times
(and $\tau^{-1}$ composed with itself $k$ times, respectively).
Let us fix a dyadic scale $\lambda=2^{j_0}$ for some  $j_0\in\Z$
and set $f_\ld(x) = f(\ld x)$.
First, note that
\begin{align*}
\big(\psi_{j, k}(D)f_\ld\big)(x)&=\int_{\xi\in 2^j(Q_0+k)}\psi(2^{-j}\xi-k)\ld^{-d}\ft f (\ld^{-1}\xi)e^{2\pi ix\cdot\xi}\,d\xi\\
&=\int_{\xi\in 2^{j-j_0}(Q_0+k)}
\psi(2^{-j+j_0}\xi-k)\ft{f}(\xi)e^{2\pi i\ld x\cdot\xi}\,d\xi\\
&=\big(\psi_{j-j_0, k}(D)f\big)(\ld x).
\end{align*}

\noi
Hence, we have
\[\|\psi_{j, k}(D)f_\ld\|_{L^p_x}
=\ld^{-\frac d p}\|\psi_{j-j_0, k}(D)f\|_{L^p_x}=2^{- j_0\frac dp}\|\psi_{j-j_0, k}(D)f\|_{L^p_x}.\]

\noi
Therefore, we obtain
\begin{align*}
\|f_\ld\|_{\M^{p, q, r}_{ \w}}&=\big\|\|f_\ld\|_{M^{p,q}_{[j]}}\big\|_{\ell^r_j(\w)}\\
&=\Big\|\big\|\|\psi_{j, k}(D) f_\ld
\|_{L_x^p(\R^d)} \big\|_{\l^q_k(\Z^d)}\Big\|_{\l^r_j(\w)}\\
&=2^{- j_0\frac d p}\Big\|\big\|\|\psi_{j-j_0, k}(D)f
\|_{L_x^p(\R^d)} \big\|_{\l^q_k(\Z^d)}\Big\|_{\l^r_j(\w)}\\
&=2^{- j_0\frac d p}
\Big\|\big\|\|\psi_{j, k}(D)f \|_{L_x^p(\R^d)} \big\|_{\l^q_k(\Z^d)}\Big\|_{\l^r_{j}(\tau^{j_0} \w)}\\
&=\ld^{-\frac dp}\|f\|_{\M^{p, q, r}_{\tau^{j_0} \w}}.
\end{align*}

\noi
Note that the vector weight on the right-hand side is now given by $\tau^{j_0} \w$.
Namely, in the general case, the scaling has an effect
of translating the weight by $\log_2 \ld$.

Let us now assume further that the vector weight ${\w}$ is multiplicative; that is, 
$w_{i+j}=w_iw_j$ for all $i, j\in\Z$, or equivalently that $w_j=w_1^j, j\in\Z$.
Under this extra assumption, we claim that the new modulation spaces $\M^{p,q, r}_{\w}$ in
Definition~\ref{DEF:mod_sc} enjoy a scaling symmetry.
Indeed, by repeating the computation above, we obtain
\begin{align}
\|f_\ld\|_{\M^{p, q, r}_{ \w}}
&=2^{- j_0\frac d p}\Big\|\big\|\|\psi_{j-j_0, k}(D)f
\|_{L_x^p(\R^d)} \big\|_{\l^q_k(\Z^d)}\Big\|_{\l^r_j(\w)} \notag \\
&=2^{- j_0\frac d p}w_{j_0}
\Big\|\big\|\|\psi_{j, k}(D)f \|_{L_x^p(\R^d)} \big\|_{\l^q_k(\Z^d)}\Big\|_{\l^r_{j}(\w)}\notag\\
&=\ld^{\log_2 (w_1)-\frac dp}\|f\|_{\M^{p, q, r}_{\w}}.
\label{X6}
\end{align}

\noi
This shows the scaling property of the  $\M^{p, q, r}_{\w}$-spaces,
when  the vector weight ${\w}$ is multiplicative.

\begin{remark}\rm
(i)  The basic property of the Fourier transform states
that modulations on the physical side correspond to translations on the Fourier side.
Unit-scale modulations, that is, unit-scale translations on the Fourier side, yield
the (usual) modulation spaces $M^{p, q}$ and the decomposition~\eqref{decomp1}.
On the other hand,  the new modulation spaces $\M^{p,q, r}_{\w}$ with scaling
are generated by modulations of all dyadic scales, that is, by translations
of all dyadic scales  on the Fourier side. Recalling that a wavelet basis is generated by all dyadic dilations and translations of a given (nice) function on the physical side,
it is natural that the family
$\big\{ \psi_{j, k}(\cdot) = \psi (2^{-j}\,\cdot\, -k)\big\}_{j \in \Z, k \in \Z^d}$,
which is essentially a wavelet basis but on the Fourier side,
appears in Definition \ref{DEF:mod_sc}.

\smallskip

\noi
(ii) As it is well known, the modulation spaces $M^{p, q}$ have an equivalent characterization via the short-time (or windowed) Fourier transform (STFT).
Given
a non-zero window function $\phi\in \S(\R^d)$,
we define  the STFT $V_\phi f$ of a tempered distribution $f\in \S'(\R^d)$ with respect to $\phi$ by\footnote{This is essentially the Wigner transform of $f$ and $\phi$.}
\[V_\phi f(x, \xi)=\int_{\R^d}f(y)\cj{\phi(y-x)}e^{-2\pi iy\cdot\xi}\,dy.\]

\noi
Then, we have the equivalence of norms:
$$\|f\|_{M^{p, q}}\sim_\phi
\nn f \nn_{M^{p, q}}:=
\|V_\phi f\|_{L^q_\xi L^p_x}=\big\|\|V_\phi f
\|_{L_x^p} \big\|_{L^q_\xi},$$
where the implicit constants depend on the window function $\phi$.

Given $j\in\Z$, let  $\dl_j f(x)=2^{-jd}f(2^{-j}x)$ and $\phi^j (x)=\phi(2^jx)$.
Then, a direct computation (see \eqref{Z3} below) shows that
\begin{align}
\|f\|_{M^{p, q}_{[j]}}
=2^{j\frac d{p'}}\|\dl_j f\|_{M^{p, q}}
\sim_\phi 2^{j\frac d{p'}}\nn \dl_j f\nn_{M^{p, q}}=2^{j\frac d {p'}}\|V_\phi (\dl_j f)\|_{L^q_\xi L^p_x}.
\label{A1}
\end{align}

\noi
Moreover, a straightforward calculation shows that
\begin{align}
\big(V_\phi (\dl_j f)\big)(x, \xi)=(V_{\phi^j}f)(2^{-j}x, 2^j\xi).
\label{A2}
\end{align}

\noi
Then, from \eqref{A1} and \eqref{A2}, we obtain
\[\|f\|_{M^{p, q}_{[j]}}\sim_\phi 2^{jd(\frac{1}{p'} + \frac 1p - \frac 1q)}\|V_{\phi^j}f\|_{L^q_\xi L^p_x}
=2^{j\frac{d}{q'}}\|V_{\phi^j}f\|_{L^q_\xi L^p_x}.\]

Based on these considerations,
if one prefers to use the $\nn \cdot\nn_{M^{p, q}}$-norm  via the STFT to define the new modulation spaces
$\M^{p, q, r}_{\w}$ with scaling,
the only difference will be in exchanging the vector weight $\w=\{w_j\}_{j\in\Z}$
with the vector weight $\ss=\{2^{j\frac{d}{q'}}w_j\}_{j\in\Z}$. Namely, if we define
\[\nn f\nn _{M^{p, q}_{[j]}}:=\|V_{\phi^j}f\|_{L^q_\xi L^p_x},\]

\noi
then we obtain the following equivalence of norms:
\[\|f\|_{\M^{p, q, r}_{\w}}\sim \nn f\nn_{\M^{p, q, r}_{\ss}}:=\big\|\nn f\nn_{M^{p,q}_{[j]}}\big\|_{\ell^r_j(\ss)}.\]

\smallskip

\noi
(iii)
The boundedness property of
the dilation operator:
$f(x)\mapsto f_\ld(x) = f(\ld x)$ 
on the modulation spaces was studied
in \cite{ST}.
In particular,
it was shown in \cite[Theorem 1.1]{ST} that
there exist $c_1, c_2 > 0$,  depending only on $d, p$, and $q$,  such that
\begin{align}
\ld^{c_1}\|f\|_{M^{p, q}}\les \|f_\ld\|_{M^{p, q}}\les \ld^{c_2}\|f\|_{M^{p, q}}
\label{A3}
\end{align}

\noi
for all $\ld \geq 1$.  When $0< \ld < 1$,
the estimate \eqref{A3} holds after switching the exponents $c_1$ and $c_2$.
The point of our construction is that the modulation spaces $\M^{p, q, r}_{\w}$
are not biased toward any particular scale.
Moreover, when the vector weight is multiplicative,
the $\M^{p, q, r}_{\w}$-spaces enjoy
the \textit{exact} dyadic scaling  \eqref{X6}, which is relevant, for example, in the analysis of the nonlinear Schr\"odinger equation; see also Sections \ref{SEC:L^2} and \ref{SEC:Application}.

\smallskip

\noi
(iv) There are variants of the estimate \eqref{A3} in the settings of weighted modulation spaces \cite{CO} or $\alpha$-modulation spaces \cite{HaW} which are interesting in their own right. 
Naturally, one could consider  a similar construction of scale invariant modulation spaces 
in these  contexts.
We, however,  do not pursue this issue here and leave it for  interested readers.

\smallskip

\noi
(v) Let $1\leq p_0 < \infty$.
Then, for $p>p_0>q'>1$, the vector weight $\w=\w(p, p_0)=\{w_j\}_{j \in \Z}$ with
\[
w_j=2^{jd (\frac 1p-\frac 1{p_0})}
\]

\noi
is a multiplicative weight satisfying  \eqref{X1} with $0<\eps<d\min \{\frac{1}{q'}-\frac 1{p_0}, \frac 1{p_0}-\frac 1p\}.$ Moreover, for this choice of weight, \eqref{X6} implies that, for a fixed dyadic scale $\lambda=2^{j_0}$,
\[
\|f_\ld\|_{\M^{p, q, r}_{\w}}=\ld^{-\frac d{p_0}}\|f\|_{\M^{p, q, r}_{\w}}.
\]

\noi
Namely,  $\M^{p, q, r}_{\w(p, p_0)}$ scales like $L^{p_0}(\R^d)$. See also Section \ref{SEC:L^2}.

\end{remark}

\section{A ``good" weight and $L^2$-embedding}\label{SEC:L^2}

In the discussion below,  we restrict our attention to the case $q=r$ and we simply write
$\M^{p, q, q}_{\w}$ as $\M^{p, q}_{\w}$.
We will show that, for a particular good vector weight $\w$ and $p, q>2$, the space $\M^{p, q}_{\w}$ is sufficiently large.

Given $ p > 2$, fix a vector weight $\w(p)  = \{w_j\}_{j \in \Z}$ with
\begin{align}
 w_j=2^{jd\frac{p' - 2}{2p'}}.
 \label{Y1}
\end{align}

\noi
We see that  $\w(p)$ is  a multiplicative vector weight and is a ``good" weight
in the sense of the conditions \eqref{X1}, provided that $q > 2$.

Before proceeding further, let us recall
 the space $X_{p,q}$  defined in
\cite[p.\,5260]{BV} via  the norm:
\begin{align}
 \|  f \|_{X_{p, q}}
=  \bigg(\sum_{j\in\Z} 2^{jd \frac{p-2}{2p}q}\sum_{k\in\Z^d}\|\psi_{j, k}f\|_{L^{p}}^q\bigg)^{\frac 1q}.
\label{Y2}
\end{align}

\noi
The $X_{p, q}$-spaces appear in the improvement of  the Strichartz estimates
for the linear Schr\"odinger equation (see \eqref{Strichartz} below) and, in particular, play an important role in the study of the mass-critical nonlinear Schr\"odinger equations.
See also \cite[Proposition 2.1]{CK} and \cite[Theorem 4.2]{MVV}
for the one-dimensional and two-dimensional versions of the $X_{p, q}$-spaces.
Note that the vector weight $\big\{ 2^{jd \frac{p-2}{2p}}\big\}_{j \in \Z}$ appearing in \eqref{Y2} is precisely  $\w(p')$ defined in~\eqref{Y1}.

Let us now come back to the $\M^{p, q}_{\w(p)}$-space
with $\w(p) = \{w_j\}_{j \in \Z}$ defined in \eqref{Y1}.
For $p > 2$, it follows from  Hausdorff-Young's inequality that
\begin{align}
\| f\|_{\M^{p, q}_{\w(p)}}
& =  \bigg(\sum_{j\in\Z}w_j^q
\sum_{ k\in\Z^d}
\|\psi_{j, k}(D)f\|_{L^p_x}^q\bigg)^{\frac 1q}\notag\\
& \le \bigg(\sum_{j\in\Z} 2^{jd \frac{p'-2}{2p'}q}\sum_{k\in\Z^d}\|\psi_{j, k}\ft{f}\|_{L^{p'}}^q\bigg)^{\frac 1q}
= \| \ft f \|_{X_{p', q}} \notag\\
& = : \| f\|_{\F X_{p', q}}.
\label{Y2a}
\end{align}

\noi
This shows that
\begin{align}
\F X_{p', q}\subset  \M^{p,q}_{\w(p)}.
\label{Y3}
\end{align}

\noi
In particular, it follows from  \cite[Theorem 1.3]{BV}
and Plancherel's theorem
that  for all $p, q>2$,  we have
\[\|f\|_{\M^{p, q}_{\w(p)}}\les \|f\|_{L^2}\]

\noi
 for all $f\in L^2(\R^d)$.
Namely, we have the following embedding:
\[L^2(\R^d)\subset \M^{p, q}_{\w(p)}(\R^d).\]

\noi
Moreover, this embedding is strict; indeed, if we set $f=\F^{-1} g$, where
\[
g(\xi)=|\xi|^{-\frac d2}\, \big|\ln|\xi| \big|^{-\frac 12}\cdot \ind_{(0, \frac 12)^d}(\xi),
\]
we see that $f\in \M^{p, q}_{\w(p)}(\R^d)\setminus L^2(\R^d)$
when $p, q > 2$; see again \cite[p.\,5267]{BV}.

\begin{remark}\rm

In the study of the one-dimensional cubic nonlinear Schr\"odinger equation
in almost critical spaces,
the \textit{Fourier-amalgam  spaces} $\ft w^{p, q} (\R^d)$,
defined by the norm:
\begin{equation*} 
\|f\|_{\ft w^{p, q} (\R^d)} = \big\|\|\psi(\xi-k) \ft f
\|_{L_\xi^p(\R^d)} \big\|_{\l^q_k(\Z^d)},
\end{equation*}

\noi
were considered in \cite{FO}. 
We note that the Fourier-amalgam space $\ft w^{p, q} (\R^d)$ is simply the Fourier image of the usual Wiener amalgam space $W(L^p, \l^q)(\R^d)$,  as described in \cite{Fei1983, FS}, defined by
the norm:
$$\|f\|_{W(L^p, \l^q)(\R^d)}=\big\|\|\psi(x-k) f
\|_{L_x^p(\R^d)} \big\|_{\l^q_k(\Z^d)}.$$

Comparing this with the definition \eqref{mod} of the modulation spaces $M^{p, q}$,
we see that $\ft w^{p', q}$ ``scales like'' $M^{p, q}$ (in the sense that
$L^{p'}_\xi(\R^d)$ scales in the same manner as  $L^p_x(\R^d)$). In Definition \ref{DEF:mod_sc}, we defined the new modulation spaces $\M^{p, q, r}_{\w}$ with scaling
for a good vector weight $\w$ satisfying \eqref{X1}.
In an analogous manner, we can also define
the new Fourier-amalgam spaces $\ft \ww^{p, q, r}_{\w}$
adapted to scaling
by the norm:
\begin{equation*}
\|f\|_{\ft \ww^{p,q, r}_{\w}}:=\big\|\|f\|_{\ft w^{p,q}_{[j]}}\big\|_{\ell^r_j(\w)}
=\bigg(\sum_{j\in\Z}
w_j^r\|f\|_{\ft w^{p,q}_{[j]}}^r\bigg)^{\frac 1r},
\end{equation*}

\noi
where the scaled Fourier-amalgam spaces $\ft w^{p,q}_{[j]}$ are defined by
\[\|f\|_{\ft w^{p,q}_{[j]}}:=\big\|\|\psi_{j, k} \ft f
\|_{L_\xi^p(\R^d)} \big\|_{\l^q_k(\Z^d)}.\]

\noi
As in the case of the $\M^{p,q, r}_{\w}$-spaces, we need to impose
some conditions on the vector weight $\w$ so that the resulting space
$\ft \ww^{p,q, r}_{\w}$ contains smooth functions.
Recalling the scaling property of $\ft m^{p, q}$,
we impose the following condition on the vector weight $\w =\{w_j\}_{j \in \Z}$:
\begin{align}
w_j\les
\begin{cases}
2^{-\eps j} & \text{if}\,\,j\geq 0,\\
2^{-   (\frac{d}{p} - \frac dq - \eps) j } & \text{if}\,\,j<0.
\end{cases}
\label{Fmod3}
\end{align}

\noi
Namely, we replace $p'$ in \eqref{X1} by $p$.
Arguing as in Subsection \ref{SUBSEC:scaling},
we see that this new Fourier-amalgam space
 $\ft \ww^{p,q, r}_{\w}$ enjoy a scaling property
analogous to that for the $\M^{p,q, r}_{\w}$-spaces.
In particular, if the vector weight $\w$ is multiplicative,
then we have the following scaling property;
if we let  $\lambda=2^{j_0}$ for some  $j_0\in\Z$ and set $f_\ld(x) = f(\ld x)$ as before, then
\begin{align*}
\|f_\ld\|_{\ft \ww^{p, q, r}_{ \w}}
&=\ld^{\log_2 (w_1)-\frac d{p'}}\|f\|_{\ft \ww^{p, q, r}_{\w}}.
\end{align*}

We conclude this remark by pointing out that the $\F X_{p, q}$-space discussed above
is simply $\ft \ww^{p, q, q}_{\w(p')}$ with the vector weight $\w(p') = \big\{ 2^{jd \frac{p-2}{2p}}\big\}_{j \in \Z}$.
It is easy to see that this vector weight satisfies the condition \eqref{Fmod3} for
$q > 2$
and it is also multiplicative, thus satisfying the following scaling property
for $\ld=2^{j_0}$ for some  $j_0\in\Z$:
\begin{align*}
\| f_\ld \|_{\F X_{p, q}}
& = \|f_\ld\|_{\ft \ww^{p, q, q}_{ \w(p')}}
=\ld^{d\frac{p-2}{2p}-\frac d{p'}}\|f\|_{\ft \ww^{p, q, q}_{\w(p')}}\\
& = \ld^{-\frac d 2}\|f\|_{\ft \ww^{p, q, q}_{\w}}
= \ld^{-\frac d 2}\|f\|_{\F X_{p, q}}.
\end{align*}

\noi
Namely, the $\F X_{p, q}$-spaces (and hence the $X_{p, q}$-spaces) scale like $L^2(\R^d)$
as it is expected.

\end{remark}

\begin{remark}\rm
In studying the mass-subcritical 
generalized KdV equation, Masaki-Segata \cite[Definition 2.1]{MS} 
 considered the so-called \textit{generalized Morrey spaces} $M^p_{q, r}$ (with $q\leq p$) 
 and $\ft{M}^{p}_{q, r}=\mathcal FM^{p'}_{q', r}$ (with $p\leq q$
 and $r > p'$).
 The usual Morrey spaces $M^p_q$ correspond to $r=\infty$ in the scale of generalized Morrey spaces, that is $M^p_q=M^p_{q, \infty}$. 
 These spaces generalize
 $X_{p, q}$
 and $\F X_{p', q}$ defined in \eqref{Y2} and \eqref{Y2a}
 and in particular we have 
 $X_{p, q} = M^2_{p, q}$ and  $\F X_{p', q} = \ft M^2_{p, q}$.
We point out that, when $r < \infty$,  the space $\ft M^p_{q, r}$
is nothing but  
 $\ft \ww^{q', r, r}_{\w(p, q)}$ with a specific vector weight given by $\w(p, q) = \big\{ 2^{jd (\frac 1q - \frac 1p)}\big\}_{j \in \Z}$,\footnote{Under $r > p'$, 
 this vector weight satisfies the condition \eqref{Fmod3} 
 for $q < 2$.}
 while $M^p_{q, r}$ is the Fourier image of 
 $\ft \ww^{q, r, r}_{\w'(p, q)}$ with the vector weight $\w'(p, q) = \big\{ 2^{jd (\frac 1p - \frac 1q)}\big\}_{j \in \Z}$.
 For more on Morrey-type spaces, 
 see also \cite{Fofana1, Fofana2, Feuto}.
\end{remark}

\section{An application to the Schr\"odinger multiplier}\label{SEC:Application}

Let  $S(t)=e^{-it\Delta}$ be the Schr\"odinger operator
defined as a Fourier multiplier operator
with a multiplier $e^{4\pi^2 it|\xi|^2}$.
Our main interest here is to discuss
a boundedness property
of the Schr\"odinger operator
 $S(t)=e^{-it\Delta}$
 on  appropriate $\M^{p, q, r}_{\w}$-spaces.

Let $\psi$ be as in \eqref{psi0}.
A direct computation shows that
\begin{align}
\psi_{j, k}(D)(S(t)f)(x)&=\int_{\R^d}\psi (2^{-j}\xi-k)e^{4\pi^2 it|\xi|^2}\ft f(\xi)e^{2\pi ix\cdot \xi}\,d\xi\notag \\
&=2^{jd}\int_{\R^d}\psi(\xi-k)e^{4\pi^2 i (2^{2j}t)|\xi|^2}\ft f(2^j\xi)e^{2\pi i2^jx\cdot\xi}\,d\xi \notag \\
&=2^{jd}\big(\psi(D-k)(S(2^{2j}t)\dl_j f)\big)(2^jx),
\label{Z0}
\end{align}

\noi
where, as before, we wrote $\dl_j$ for the dilation operator given by $\dl_j f(x)=2^{-jd}f(2^{-j}x)$.
Thus, we have
\[
\|\psi_{j, k}(D)(S(t)f)\|_{L^p_x }=2^{j\frac{d}{p'}}\|\psi(D-k)(S(2^{2j}t)\dl_j f) \|_{L^p_x }.
\]

\noi
Hence, we obtain
\begin{align}
\|S(t)f\|_{M^{p, q}_{[j]}}=2^{j\frac{d}{p'}}\|S(2^{2j}t)\dl_j f\|_{M^{p, q}}.
\label{Z1}
\end{align}

\noi
Interpolating   \cite[Theorem 14]{BGOR} for $p = 1$ or $p = \infty$
with the trivial bound for $p = 2$, we obtain
\begin{align}
\|S(2^{2j}t)\dl_j f\|_{M^{p, q}}\les \jb{2^{2j}t}^{d\,|\frac 12-\frac 1p|}\|\dl_j f\|_{M^{p, q}}.
\label{Z2}
\end{align}

\noi
Lastly, by a computation analogous to \eqref{Z0}, we have
\begin{align*}
\psi(D-k)(\dl_j f)(x)
= 2^{-jd}\psi_{j, k}(D)f(2^{-j}x)
\end{align*}

\noi
and hence
\begin{align}
\|\dl_j f\|_{M^{p, q}}
=2^{-jd}2^{j\frac dp}\big\|\|\psi_{j, k}(D)f\|_{L^p_x}\big\|_{\l^q_k}
=2^{-j\frac{d}{p'}}\|f\|_{M^{p, q}_{[j]}}.
\label{Z3}
\end{align}

\noi
Putting \eqref{Z1}, \eqref{Z2}, and \eqref{Z3} together, we obtain
\begin{align}
\|S(t)f\|_{M^{p, q}_{[j]}}\lesssim \jb{2^{2j}t}^{d\, |\frac12-\frac1p|}\|f\|_{M^{p, q}_{[j]}}.
\label{Z4}
\end{align}

\noi
Note that, when $j=0$, this estimate recovers the boundedness of $S(t)$ on the modulation spaces $M^{p, q}$.
We also point out that, when $p \ne 2$,
the divergent behavior of the constant as  $j \to \infty$
is consistent with the unboundedness of $S(t)$ on $L^p(\R^d)$.

Now, let  $\w = \{w_j\}_{j \in \Z}$ be a ``good" vector weight, satisfying \eqref{X1},
and define  a vector weight $\ss=\{\s_j\}_{j\in\Z}$ by setting
\[ \s_j = \begin{cases}
2^{- j d\, | 1 - \frac 2p|}w_j & \text{if } j \geq 0,\\
w_j & \text{if } j < 0.
\end{cases}
\]

\noi
Note that such $\ss$ also satisfies \eqref{X1}.
Then,
from \eqref{Z4},  we obtain the boundedness of
$S(t): \M^{p, q, r}_{\w}\to \M^{p, q, r}_{\ss}$ with a quantitative bound:
\[
\|S(t)f\|_{\M^{p, q, r}_{\ss}}\les \jb{t}^{d\, |\frac 12 - \frac 1p|}\|f\|_{\M^{p, q, r}_{\w}}
\]

\noi
for any $t \in \R$.

Let us further comment on an alternate estimate when $p\geq 2$. In this case, the Schr\"odinger operator is known to also satisfy the following inequality (see \cite[Proposition 4.1]{HW}):
\begin{equation*}
\|S(t)f\|_{M^{p, q}}\lesssim \jb{t}^{-d\,(\frac 12 - \frac 1p)}\|f\|_{M^{p', q}}.
\end{equation*}

\noi
Note that the decay at infinity in this estimate is now consistent with the one in the $L^{p'}\to L^p$ boundedness of $S(t)$.
Thus, instead of \eqref{Z2}, we can now write
\begin{align}
\|S(2^{2j}t)\dl_j f\|_{M^{p, q}}\les \jb{2^{2j}t}^{d\,(\frac 12-\frac 1p)}\|\dl_j f\|_{M^{p', q}}.
\label{Z5}
\end{align}
Since, by \eqref{Z3}, 
we have $\|\dl_j f\|_{M^{p', q}}=2^{-j\frac{d}{p}}\|f\|_{M^{p', q}_{[j]}},$ \eqref{Z1} and \eqref{Z5} yield
\begin{align}
\|S(t)f\|_{M^{p, q}_{[j]}}&\lesssim \jb{2^{2j}t}^{-d\, (\frac12-\frac1p)}2^{jd\, (\frac1{p'} - \frac 1p)}\|f\|_{M^{p', q}_{[j]}}\notag\\
&=\jb{2^{2j}t}^{-d\, (\frac12-\frac1p)}2^{2jd\, (\frac12 - \frac 1p)}\|f\|_{M^{p', q}_{[j]}}.
\label{Z6}
\end{align}
Thus, by considering again $\w = \{w_j\}_{j \in \Z}$ to be a ``good" vector weight, satisfying \eqref{X1}, from \eqref{Z6} we obtain now the boundedness $S(t): \M^{p', q, r}_{\w}\to \M^{p, q, r}_{\w}$ with a quantitative bound:
\[
\|S(t)f\|_{\M^{p, q, r}_{\w}}\les \jb{t}^{-d\, (\frac 12 - \frac 1p)}\|f\|_{\M^{p', q, r}_{\w}}
\]

\noi
for any $p\geq 2$ and $t \in \R$.

We conclude our note with a comment on the following improvement
of the Strichartz estimate proved in \cite[Theorem 1.2]{BV}; if $2<p<2+\frac{4}{d(d+3)}$ and $q=\frac{2(d+2)}{d}$, then
\begin{equation}
\label{Strichartz}
\|S(t)f\|_{L^q_{t, x} (\R\times \R^{d})}\lesssim \|f\|_{\F X_{p', q}(\R^{d})}.
\end{equation}

\noi
The proof of this inequality is highly non-trivial; it uses the bilinear restriction estimate in \cite[Theorem 1.1]{Tao} and an appropriate orthogonality lemma for functions with disjoint Fourier supports \cite[Lemma 6.1]{TVV}.
Considering the embedding \eqref{Y3} stated in the previous section, it would be interesting to know
if the $\F X_{p', q}$-norm in \eqref{Strichartz} can be replaced by  the $\M_{\w(p)}^{p, q}$-norm.

\begin{ackno}\rm
The first author is partially supported by a grant from the Simons Foundation (No.\,246024). The second author is supported by the European Research Council (grant no.~637995 ``ProbDynDispEq'').
The authors would like to thank the anonymous referees for several useful comments
that have improved the presentation of this paper.

\end{ackno}

\end{document}